\documentclass[12pt]{article}
 \usepackage{amsmath,amsfonts,theorem}
 \numberwithin{equation}{section}

 \newtheorem{prop}{Proposition}[section]

 \newtheorem{cor}{Corollary}[section]

 \newtheorem{thm}{Theorem}[section]

 \theorembodyfont{\rmfamily}
 \newtheorem{dfn}{Definition}[section]

 \newtheorem{rmk}{Remark}

 \newcommand{\qed}{\ifhmode\unskip\nobreak\fi\quad\ensuremath\square}

 \newcommand{\p}{\partial}

 \newcommand{\ov}{\overline}

 \newcommand{\sA}{\mathcal A} 
 
 \newcommand{\sG}{\mathcal G} 
 \newcommand{\sC}{\mathcal C}

 \newcommand{\sM}{\mathcal M}
 \newcommand{\Oh}{\mathcal O}
 \newcommand{\sP}{\mathcal P}
 \newcommand{\sS}{\mathcal S}
 
 \newcommand{\ep}{\varepsilon}
 
 \newcommand{\ga}{\gamma}
 
 \newcommand{\om}{\omega}
 \newcommand{\si}{\sigma}
 
 \newcommand{\Ga}{\Gamma}
 
 \newcommand{\Om}{\Omega}
 \newcommand{\la}{\lambda}
 \newcommand{\Si}{\Sigma}

 \newcommand{\PP}{\mathbb P}
 \newcommand{\C}{\mathbb C}

 \newcommand{\Z}{\mathbb Z}

 \newcommand{\CLRep}{\operatorname{CLRep}}
 \newcommand{\SL}{\operatorname{SL}} 

 \newcommand{\SU}{\operatorname{SU}}

\begin{document}

 \title{On the Mumford-Narasimhan problem}
 \markright{\hfill Large limits... \quad}

 \author{Andrei Tyurin}
 \maketitle

\begin{center}
 {\em To Prof. M. Narasimhan on his 70th birthday}
 \end{center}
 \bigskip

Thirty years ago, after appearing of the papers \cite{T1} David Mumford
asked me about an application of the proposed  technique to the 
"Schottki problem for vector bundles". This problem is a non-abelian analogy
 of the classical cover $f \colon \C^g \to J(C)$ where the target space is the place where the theta functions of an algebraic curve $C$ live. More precisely, if $S^r_g  \subset \CLRep (\pi_1(C), \SL(r, \C))$ is the subset of the representations with trivial  "a"-periods and $f \colon S^r_g \to M_r^{ss}(C)$ is the forgetful map from the space of flat bundles to the moduli space of semi-stable vector bundles  on $C$ then the question is what kind of image we have. The spectrum of the cases is the following
\begin{enumerate}
 \item  surjectivity $f(S^r_g) =  M_r^{ss}(C)$ is the best answer, predicted by the classical case;
\item for  general curve $f(S^r_g)$ contains a Zariski open set in $ M_r^{ss}(C)$;
\item other possibilities are too bad to discuss them.
\end{enumerate}
Many times after we discussed this problem with M. Narasimhan
reducing the question to new approaches like Hitchin's Higgs bundles \cite{H} and so on. The point was the  following:
 by the Narasimhan-Sesadri theorem the restriction of $f$ to the unitary Schottki space $$uS^r_g = \CLRep(\pi_1(C), \SU(r)) \subset S^r_g$$ is an embedding. Moreover 
 the differential of $f$ can be investigated. C. Florentino proved that around $uS^r_g$ this differential is an isomorphism. Long time after we hadn't any information about its behavior. In this paper we prove at least that the second case is realized. That is  for  general smooth curve $C$  general vector bundle admits the Schotki representation. It is quite relevant to dedicate this paper to Prof.
M.Narasimhan on his 70th birthday in spite of the fact that this is the first step only to an expected perfect  solution to  the Mumford-Narasimhan problem. The main idea is 
\begin{enumerate}
\item to prove the statement for special curves and
\item  to extend
 this statement using the standard technique of complex analysis to  neighborhoods of special curves.
\end{enumerate}
 This problem lies on the boundary connecting algebraic geometry and complex analysis. But  for our special curves the  analysis can be reduce to a
simple algebraic geometry. These curves are called {\it large limit}
 curves and we would not discuss here why. But it is quite easy to understand from section 3 which isn't necessery for the main result proof. Such curves were used succesfully by Ciro Ciliberto,
Angelo Lopez, Rick Miranda and Lucia
Caporaso  already in algebro-geometrical set up (see, for example \cite{CLM}). Here we use them in  the complex gauge theoretical set up. Moreover,  the correlation functions of all local quantum
  fields in a two dimensional conformal field theory can be recovered from the partition function when all
  channels (tubes of the pair of pants decomposition)  of the surface are constricted down to nodes. This procedure  produces our reducuble curve $P_{\Ga}$ with uniquely defined complex structure. Details of such approach can be found in \cite{T3}.

\section{ Trinion decompositions and holomorphic flat connections on smooth algebraic curves}

For an oriented compact surface $\Si$ of genus $g$ 
a trinion decomposition is given by  a maximal set of 
disjoint, isotopy  inequivalent circles 
\begin{equation}
( C_1, ... , C_{3g-3} ) \subset \Si
\end{equation}
(see \cite{HT}). Removing these circles we get
\begin{equation}
\Si - \{ C_1, ... , C_{3g-3} \} = \cup_{i=1}^{2g-2} \tilde{v_i}
\end{equation}
  the finite set of trinions (or "pairs of pants") that is a trinion decomposition of $\Si$. Such decomposition 
defines   {\it dual
 trivalent graph} $\Ga$ such that
\begin{enumerate}
\item  the set of its vertices 
\begin{equation}
V(\Ga) = \{ v_i \} = \{ \tilde{v_i} \}
\end{equation}
is the set of trinions;
\item  and with the set of edges
\begin{equation}
E(\Ga) = \{ e \} = \{ C_e \};
\end{equation}
\item and two vertices $v$ and $v'$ are joined by 
 edge $e$ iff there exists the  corresponding circle
\begin{equation}
C_e = \p \tilde{v} \cap \p \tilde{v'}.
\end{equation}
\end{enumerate}

Actualy we can start with any  3-valent graph
 $\Ga$  with the set of edges $E(\Ga)$, then $|E(\Ga)| = 3g-3$, and with the set of vertices 
$V(\Ga), \quad |V(|\Ga)| = 2g-2$ where $|S|$ is the cardinality of a finite set $S$.

Topologically it is equivalent to a
3-dimensional handlebody $H_{\Ga}$ with boundary
$\p H_{\Ga} = \Si_{\Ga}$ 
where $\Si_{\Ga}$ is the Riemann surface given by the pumping up
trick (see \cite{T2}): we pump up every edge of $\Ga$ to a tube and every vertex to a trinion (=  2-sphere with 3
holes). By the construction our Riemann surface $\Si_{\Ga}$ has a
trinions (or  "pair of pants") decomposition given by removing all
tubes.

For a visualization we can input our graph 
\begin{equation}
i \colon \Ga \hookrightarrow H_\Ga
\end{equation}
in the handlebody by the natural way. Then we can see that \begin{enumerate}
\item an orientation $\vec e$ of an edge $e \in E(\Ga)$ gives the orientation $\vec C_e$ of the corresponding circle from the collection     (1.1);  
\item sending any 1-cycle on $\Si_\Ga$ to the cycle on the handlebody $H_\Ga$ we get an epimorphism
 of the fundamental groups:
\begin{equation}
r \colon  \pi_1 (\Si_\Ga) \to \pi_1 (\Ga) = \pi_1(H_\Ga) \to 1
\end{equation}
 the kernel of which is the  free group with $g$ generators
\begin{equation} 
ker \quad r = F_g.
\end{equation}
\item  According to the corresponding exact sequence
\begin{equation}
1 \to ker \quad r \to \pi_1 (\Si_\Ga) \to \pi_1 (\Ga) = \pi_1(H_\Ga) \to 1
\end{equation}
 we can choose  standard generators of $\pi_1 (\Si_\Ga)$ 
$$
\pi_1(\Si_\Ga) = <a_1, ... ,a_g, b_1, ... , b_g \vert \prod_{i=1}^g [a_i, b_i] = 1>
$$
such that
\begin{equation}
ker \quad r = F_g = <a_1, ... , a_g>
\end{equation} and
$$
\pi_1 (\Ga) =  <r(b_1), ... , r(b_g)>;
$$
\item obviously whole the collection of cycles $\{ [C_e]\}$     (1.1) lies in $ker \quad r$ and moreover
\begin{equation}
ker \quad r = F_g = < [C_1], ... , [C_{3g-3}]> ;
\end{equation}
\end{enumerate}

Recall that a path of lenth 1 on  $\Ga$ is just an oriented edge
 $\vec e$. Let the set
$P_1(\Ga) = \vec E(\Ga)$ be the set  of  1-paths on  $\Ga$.
Every such path  $\vec e $ has  vertices of two types
\begin{equation}
v_s(\vec e), \quad v_t(\vec e) \in V(\Ga)
\end{equation}
- the {\it source} and the {\it target} - which are equal for  loops.

A path of length  $d$ in  $\Ga$ is an ordered sequence $(\vec e_1,
..., \vec e_d)$ of oriented edges such that for every $i$
\begin{equation}
v_t(\vec e_i) = v_s(\vec e_{i+1}).
\end{equation}
If $\vec e_{d+1} = \vec e_1$, then our path is a loop. A path $(\vec
e_1, ..., \vec e_d)$ (or a loop)    is called irreducible if   $e_i \neq
e_{i+1}  $
 for every
$i$ (including $i = d+1 $ for a loop).

Every path
 $(\vec e_1, ..., \vec e_d) \in P_d(\Ga)$ defines two vertices
\begin{equation}
v_s((\vec e_1, ..., \vec e_d)), \quad v_t((\vec e_1, ..., \vec
e_d)) \in V(\Ga)
\end{equation}
- the source and  the target and the coincidence of which shows that our parth is a loop.

Let  
\begin{equation}
L_d(\Ga) \subset  P_d(\Ga)
\end{equation}
be the set of oriented irreducible loops of length   $d$.  For every
vertex á
 $v \in V(\Ga)$ we have the set of irreducible loops marked by $v$
\begin{equation}
L_d(\Ga)_v = \{ l \in L_d(\Ga) \vert v \subset l\}.
\end{equation}
The union
\begin{equation}
L_\infty(\Ga)_v = \cup_{d=1}^\infty L_d(\Ga)_v
\end{equation}
admits a group structure
\begin{equation}
\pi_1^C(\Ga)_v = L_\infty(\Ga)_v.
\end{equation}
if we consider only irreducible fragments of the compositions.

Obviously, this group depends on the marking point $v$.

Let $\pi_1(\Ga)$ be the standard fundamental group of $\Ga$ (as a
1-complex). Then the natural epimorphism  $r \colon \pi_1^C(\Ga)
\to \pi_1(\Ga)$ is an isomorphism. Obviously, if $\Ga$ is a 3-valent graph 
of genus $g$ then $\pi_1(\Ga) = F_g $
is  a free group with $g$ generators.

Now let us fix a point $p \in \Si_\Ga$ (suppose it coincides with $i(v)$     (1.6)) and join all circles with $p$ by any system of paths. Then we can consider the oriented
circles $\vec C_e$  as elements of the group $ker \quad r$     (1.9) that is as elements of the fundamental group $\pi_1(\Si_\Ga)$. 

Sending every combinatorial loop 
$$
l = (\vec e_1, ..., \vec e_d) \in  L_\infty(\Ga)_v
$$
to the ordered product
\begin{equation}
Int(l) = \vec C_{e_1} \cdot ... \cdot \vec C_{e_d} \in ker \quad r \subset \pi_1(\Si_\Ga)
\end{equation}
we get a homomorphism
\begin{equation}
Int \colon \pi_1(\Ga) = \pi_1(H_\Ga) \to \pi_1(\Ga)
\end{equation}
Recall that both of these groups are free groups with $g$ generators and it is easy to see that $Int$ is an isomorphism.

Our surface $\Si_\Ga$ doesn't carry any "natural" complex structure. Moreover after the fixing of the previous basis of the fundamental group     (1.9)-    (1.10) the space of complex structures on
$\Si_\Ga$ turns to be the Teichmuller space $\tau_g$. If we fix 
such a structure $I \in \tau_g$ then we can consider the space of classes of
representations 
\begin{equation}
A^{na} =  CLRep(\pi_1(\Si), SL(2, \C))
\end{equation}
 as the full space $A^{na}(\Si_I)$ of all
 holomorphic flat connections on  topologicaly trivial
vector bundle on $\Si_I$. Any representation $\rho$
 with  class
\begin{equation}
[\rho] \in  CLRep(\pi_1(\Si), SL(2, \C))
\end{equation}
defines a holomorphic vector bundle
just by the standard construction
\begin{equation}
E  = U \times \C^2 / (\pi_1,\rho)
\end{equation}
 where $U$ is the universal cover of $\Si_I$ with the natural
  action of the fundamental group of the base.

Of course, in the set of such bundles there are non-stable but
indecomposable
 vector bundles and even semi-stable  decomposable ones. Let us remove
the representations corresponding to non-stable bundles and get the space 
\begin{equation}
 CLRep^{I-ss}(\pi_1(\Si_I), \SL(2, \C)) \subset  CLRep(\pi_1(\Si), \SL(2, \C))
\end{equation}
of classes
 of representations which give stable, semi-stable and  decomposable  semi-stable bundles.

 Obviously this Zariski dense (=containing a Zariski open set) subspace depends on the complex structure $I \in \tau_g$. Moreover we have the       holomorphic forgetful map
\begin{equation}
 f \colon A^{na}_{I-ss} = CLRep^{I-ss}(\pi_1(\Si),  SL(2, \C)) \to M^{ss}(\Si_I)
\end{equation}
where the target space is the moduli space of semi-stable 
topologicaly trivial vector bundles.
This is the affine bundle over the cotangent bundle of $ M^{ss}(\Si_I)$. Indeed, any fiber 
$f^{-1}(E)$ is the affine space of holomorphic flat connections on $E$ and the difference of any two connections is
a traceless Higgs field
\begin{equation}
\phi \colon E \to E (K_{\Si_I}) \in T^*_{\Si_I} M^{ss}(\Si_I)
\end{equation}
that is a  covector to $M^{ss}(\Si_I)$ at $E$.

 Every affine  bundle is given by a 1-cocycle of the vector bundle (see \cite{T1}). In our case this is a cocycle 
\begin{equation}
\ep_{na} \in H^1(M^{ss}, \Om) 
\end{equation}
and 
$$
 H^1(M^{ss}, \Om) = H^{1,1}(M^{ss}, \C)
$$
by the Dolbault theorem. In our case this class is precisely the class of the polarization that is the class of theta divisor:
\begin{prop} The cohomology class (1.27) is given by the formula
\begin{equation}
\ep_{na} = [\Theta_{na}] = [\om_{na}]
\end{equation}
where $\om$ is the standard symplectic structure on the moduli space.
\end{prop}

\begin{cor} The forgetful projection $f$     (1.25) doesn't admit any holomorphic section.
\end{cor}
That is, the answer is precisely the same as in the abelian case.
\begin{rmk}
By the Narasimhan-Sesadri theorem there exists the non-holomorphic section
\begin{equation}
M^{ss}(\Si_\Ga) =  CLRep(\pi_1(\Si), \SU(2)) \hookrightarrow  CLRep(\pi_1(\Si), \SL(2, \C))
\end{equation}
\end{rmk}
Now let us perform the trinion decomposition corresponding
to $\Si_\Ga$     (1.1)-    (1.2). We can construct the space of
holomorphic flat bundles $ CLRep(\pi_1(\Si), \SL(2, \C))$
 gluing the spaces of $\SL(2, \C)$-flat connections over the trinions $\{ \tilde{v_i}\}$     (1.2).

For one trinion $\tilde{v}$ (which is 2-sphere with 3 holes)  the space of $\SL(2, \C)$-flat connections is given as the classes representations space
\begin{equation}
A^{na}(\tilde{v}) =  CLRep(\pi_1(\tilde{v}), \SL(2, \C))
\end{equation}
where $\pi_1(\tilde{v}) = F_2$ is the free group with 2 generators. 

\begin{dfn} For the free group $F_g$ with $g$ generators
 the classes representations space 
\begin{equation}
S_g =   CLRep(F_g, \SL(2, \C))
\end{equation}
is called the Schottki space of genus $g$.
\end{dfn}

So $A^{na}(\tilde{v}) = S_2$ is the Schottki space of genus 2. 

If two vertices $v_i$ and $v_j$ are joined by 
 edge $e_l$ that is if the circle $C_l$ corresponding to $e_l$ is a boundary component of trinions $\tilde{v_i}$ and  $\tilde{v_j}$: 
\begin{equation}
C_l = \p \tilde{v_i} \cap \p \tilde{v_j}
\end{equation}
then we can glue $A^{na}(\tilde{v_i}) $ and $A^{na}(\tilde{v_j}) $ along the boudary component $C_l$
by just the same way as we do this for $\SU(2)$ representations (see for example \cite{JW}): let $Conj(\SL(2, \C))$ be the set of conjugation classes of elements of $\SL(2, \C)$. Then we have the natural map
\begin{equation}
conj \colon \SL(2, \C) \to Conj(\SL(2, \C)).
\end{equation}
Two representations $\rho_i, \quad [\rho_i] \in A^{na}(\tilde{v_i})$ and   $\rho_j, \quad [\rho_j] \in A^{na}(\tilde{v_j})$ are glued iff 
\begin{equation}
conj(\rho_i([C_l])) = conj(\rho_j([C_l])).
\end{equation}
The ambiguity of such gluing is the stabilizer $Z(\rho_i([C_l]))$ of the monodromy around this 
loop. For example for a semi-simple element
$m \in \SL(2, \C)$ the stabilizer $Z(m) = \C^*$.  The result of such gluing is
\begin{equation}  
A^{na}(\tilde{v_i}) * A^{na}(\tilde{v_j}) = S_3
\end{equation}
if $C_l = \p \tilde(v_i) \bigcap \p \tilde(v_j) $.
Gluing all trinions $\tilde{v_i}$ we get our surface and gluing all spaces of flat bundles $A^{na}(\tilde{v})$ we get the space $ CLRep(\pi_1(\Si), \SL(2, \C))
$ of flat bundles on $\Si_\Ga$.

The exact sequence     (1.9) defines the canonical embedding of the Schottki space      (1.31)
\begin{equation}
i \colon S_g \hookrightarrow  CLRep(\pi_1(\Si), \SL(2, \C))
\end{equation}
such a way that the image is in a sense a complete intersection. Namely, consider the space of representations $Rep(\pi_1(\Si), \SL(2, \C))$ ( without the diagonal adjoint factorization), such that
\begin{equation}
\CLRep(\pi_1(\Si), \SL(2, \C)) = Rep(\pi_1(\Si), \SL(2, \C)) / PGL(2, \C)
\end{equation}
where the last action is the  natural diagonal adjoint action of $\SL(2, \C)$ on the space of representations. 
So, let 
\begin{equation} 
p \colon Rep(\pi_1(\Si), \SL(2, \C)) \to \CLRep(\pi_1(\Si), \SL(2, \C))
\end{equation}
be the natural projection to the factor by this action.

Geometrically this means that we fixed a point $p \in \Si_I$ and  a trivialization $E_p = \C^2$ of the fiber of vector bundle $E$ over this point. The moduli space of vector bundles with such additional structure $\widetilde{M^{ss}(\Si_I)}$ admits the structure of  principal $PGL(2, \C)$-bundle over $M^{ss}(\Si_I)$: 
\begin{equation}
\phi \colon \widetilde{M^{ss}(\Si_I)}   \to M^{ss}(\Si_I)
\end{equation}
where the group $\SL(2, \C)$ modulo $\pm 1$ acts on  $E_p = \C^2$ as (2 X 2)-matrices with determinant 1.  

The moduli space of flat vector bundles with such additional structure is $Rep(\pi_1(\Si), \SL(2, \C))$ and it admits the structure $p$ (1.38) of a principal $PGL(2, \C)$-bundle over $\CLRep(\pi_1(\Si), \SL(2, \C))$. 

Now every element $\ga \in \pi_1(\Si)$ defines  regular functions 
\begin{equation}
c_{ij}(\ga) \colon  Rep(\pi_1(\Si), \SL(2, \C)) \to \C
\end{equation}
 by the formula
$$
c_{ij}(\ga)(\rho) =  (\rho(\ga))_{ij}
$$
 where  $ij$ is the matrix element of the corresponding matrix.

Then the half basis $(a_1, ... , a_g)$     (1.10) defines $3g$ functions
\begin{equation}
c_{11}(a_1), c_{12}(a_1), c_{21}(a_1), ... , c_{11}(a_g),
c_{12}(a_g), c_{21}(a_g)  
\end{equation}
and a system of divisors 
\begin{equation}
D_{ij}(a_l) = \{ c_{ij}(a_l)  = 0 , \text{ if } i\neq j,    \}
\end{equation}
 and 
$$
D_{11}(a_l) = \{ c_{11}(a_l)  = 1 \}. 
$$
These divisors are not invariant with respect to the diagonal adjoint action of $\SL(2, \C)$. But the complete
 intersection
\begin{equation}
 \bigcap_{ij, l} D_{ij}(a_l) = \phi^{-1}(i(S_g))
\end{equation}
is the preimage of the Schottki space (1.36). 

Returning to the affine bundle     (1.25) consider a
 bundle $E \in M^{ss}(\Si_I)$ and the fiber
\begin{equation} 
f^{-1}(E) = \C^{3g-3}_E
\end{equation}
Such affine subspace of $ CLRep(\pi_1(\Si), \SL(2, \C))
$ is called a {\it packet} of classes of representations. This is an affine space over the space of Higgs fields $H^0(\Si_I, ad E \otimes T^*\Si_I)$.

Now consider the restriction of the principal $\SL(2, \C)$-bundle
(1.38) to this affine space
\begin{equation}
p \colon  p^{-1} (\C^{3g-3}_E) \to \C^{3g-3}_E.
\end{equation}
Obviously the space $ p^{-1} (\C^{3g-3}_E)$ is an affine algebraic variety. 

Restricting the functions (1.40) to this affine variety
 we get the complete intersection
\begin{equation}
 \phi^{-1}(i(S_g)) \cap p^{-1} (\C^{3g-3}_E) 
\end{equation}
and 
{\it $E$ admits a Schottki representation iff this complete intersection is non empty}.

Non emptyness is an (Zariski) open condition thus
to solve the Mumford-Narasimhan problem for  general curve  and general vector bundle  we need to find one curve with this property. 

That is our strategy is the same as C. Ciliberto, A. Lopez and R. Miranda  in \cite{CLM} for the computation of corank of the Gaussian-Wahl map for general curve. Moreover, the probe curves are the same.

\section{3-valent graphs and ll-curves}

 Sending every vertex $v \in
V(\Ga)$ to the complex Riemann sphere $P_{v} = \C\PP^{1}$ with a
triple of points $(p_{e_1}, p_{e_2}, p_{e_3})$ corresponding to edges of the star $S(v)$ of the vertex (= the set of edges incident to $v$) and gluing two components $P_{v}$ and $P_{v'}$ such that  $\p e = v, v'$ along the points $ p_{e} \in P_v, P_{v'}$ 
 we obtain a reducible algebraic curve $P_\Ga$ with the following
properties:
\begin{enumerate}
\item  arithmetical genus of the connected reducible curve
 $P_{\Ga_{g}}$ is equal to $g$;
 \item for any 3-valent graph $\Ga_{g}$ the curve $P_{\Ga_{g}}$ is
  Deligne - Mumford stable 
\item hence the curve $P_{\Ga_{g}}$ defines a point of the Deligne-Mumford
   compactification $\ov{\sM_{g}}$ (with the same notation);
\end{enumerate}

So in the Deligne-Mumford compactification $\ov{\sM_{g}}$ of the
moduli space $\sM_{g}$ of smooth curves of genus $g$ we get a finite
configuration of points
\begin{equation}
\sP \subset \ov{\sM_{g}}
\end{equation}
enumerated by the set of 3-valent
 graphs $TG_{g}$ - the set of large limit curves.

In the algebraic geometry framework these curves were described by Artamkin in \cite{A}:
\begin{enumerate} 
\item  on every curve $P_{\Ga}$ the canonical class $K_{\Ga}$ is a line
bundle for which  the  restriction  to every component $P_{v}, v \in
V(\Ga)$ is the sheaf of meromorphic differentials $\om$ with
simple poles at $p_{e_1}, p_{e_2}, p_{e_3}$ where $e_1 \bigcup e_2 \bigcup e_3 =
S(v)$:
\begin{equation}
K_{\Ga} \vert_{P_{v}} = K_{P_{v}}(p_{e_1} + p_{e_2} + p_{e_3} ) =
\Oh_{P_{\Ga}}(1);
\end{equation}
\item thus every holomorphic section $s$ of $K_{\Ga}$ is a collection of
meromorphic differentials $\{ \om_{v}\}$ on the components $P_{v}$
with poles at $\bigcup_{e \in E(\Ga)} p_{e}$ with constraints: for every $e$ such that $\p e = v, v'$
one has
\begin{equation}
res_{p_{e}}\om_{v} +  res_{p_{v'}} \om_{v'} = 0;
\end{equation}
\item beside of these equalities we have 2g-2 linear relations: for every $v \in
V(\Ga) $ with $(e_1,e_2, e_3 )= S(v)$
\begin{equation}
res_{p_{e_1}} \om_{v} + res_{p_{e_2}} \om_{v} + res_{p_{e_3}}
\om_{v} = 0;
\end{equation}
\item from here we have the right number $g$
 for  the dimension of holomorphic differential space $H^0(P_\Ga), \Oh (K_\Ga))$;
\item but the properties of the canonical complete linear system $\vert K_\Ga \vert$ depend on the topology of $\Ga$.
\end{enumerate}

\begin{dfn} The minimal number of edges that may be removed to make the graph disconnected is called  the {\it thickness} $th(\Ga)$ of the graph.
\end{dfn}
Obviously for a 3-valent graph $\Ga$ the number $th(\Ga) \leq 3$.

 In \cite{A} Artamkin proved the following

\begin{prop} The canonical linear system $\vert K_{P_\Ga} \vert$ is
\begin{enumerate}
\item base points free iff $th(\Ga) \geq 2$;
\item very ample iff  $th(\Ga) \geq 3$.
\end{enumerate}
\end{prop}

The double canonical systems on ll-curves are much more regular:
\begin{enumerate}
\item  every holomorphic section $s$ of $\Oh(2K_{\Ga})$ is a collection of
meromorphic quadratic differentials $\{ \om_{v}\}$ on the components $P_{v}$
with poles of degree 2 at $\bigcup_{e \in E(\Ga)} p_{e}$;
\item every such quadratic differential defines {\it biresidue} $bires_{p_e}$ at every pole ($Res^2_{p_e}$ in notations of \cite{T1}) and 
 \item the constraints in this case are the following: for every $e$ such that $\p e = v, v'$
\begin{equation}
bires_{p_{e}}\om_{v} +  bires_{p_{v'}} \om_{v'} = 0;
\end{equation}
\item thus the system of nodes $\{ p_e \}$ of our ll-curve $P_\Ga$ defines the decomposition
\begin{equation}
H^0(P_\Ga, \Oh (2K_\Ga))^* = \C^{E(\Ga)}
\end{equation}
where for every quadratic differential $\om$
the value of the linear form $H_e$ is given by
the formula 
\begin{equation}
H_e (\om) = bires_{p_e} \om;
\end{equation}
\item thus for a  3-valent graph $\Ga_{g}$ the large limit curve $P_{\Ga}
\in \ov{\sM_{g}}$ is a smooth point of $\ov{\sM_{g}}$ as orbifold.
\end{enumerate}

Recall that the fiber of the tangent bundle
of the moduli space at a point $P_\Ga$ is given by the equality
\begin{equation}
T_{P_{\Ga}} \ov{\sM_{g}} = H^{0}(P_{\Ga}, \Oh_{P_\Ga} (2K_{P_\Ga})
\end{equation}
and the decomposition    (2.6) is given by the following geometrical way:  the double  canonical map of $P_\Ga$ given  by the
complete linear system $\vert 2K_{\Ga}\vert$
\begin{equation}
\phi_{2K_{\Ga}} \colon P_{\Ga} \to \PP^{3g-4}
\end{equation}
has as the target space the projectivization of the tangent space of the moduli space $ \ov{\sM_{g}}$ at the point $P_\Ga$.

Then the images  of nodes 
\begin{equation}
\{ \phi_{2K_{\Ga}}(\bigcup_{e \in E(\Ga)} p_{e}) \}
\end{equation}
define the configuration of 
$$3g-3 = rk H^0(P_\Ga, \Oh(2K_{P_\Ga}))$$ 
linear independent  points ( since the components are coming to irreducible conics
and every conic is defined by any triple of points on it).  Thus this configuration of points gives a decomposition of the tangent space
\begin{equation}
T_{P_{\Ga}} \ov{\sM_{g}}= \bigoplus_{e \in E(\Ga)} \C_{e}
\end{equation}
where $\PP \C_{e} = \phi_{2K_{\Ga}}(p_{e})$
coincides with  the decomposition    (2.6).

Let $ e \subset \Ga$ be an edge of a 3-valent graph $\Ga$ with
two vertices  $v, v' = \p e$ and with two  stars $S(v)= e, e_{1}, e_{2}$ and $S(v') =
e, e'_{1}, e'_{2}$. From this we can get a graph
$\Ga'_e$ of genus $g-1$ by the following construction: 
\begin{enumerate}
\item remove $e$ and get two 2-valent vertices $v$
 and $v'$ with stars $S(v) = e_1, e_2$ and  $S(v') = e'_1, e'_2$;
\item consider the pair $ e_1, e_2$ as the first new edge
$e_{new}$ and  $e'_1, e'_2$ as the second new edge
$e'_{new}$;
\item so, we get a new graph $\Ga'$ with fixed
 pair of (a priori) disjoint edges $(e_{new} \bigcup e'_{new}) \subset \Ga'$;
\item obviously genus of $\Ga'$ is equal to $g-1$.
\end{enumerate}

Coming to ll-curves we have the correspondence:
\begin{equation}
(p_e \in P_\Ga) \longleftrightarrow (p_{e_{new}}, p_{e'_{new}} \subset P_{\Ga'})
\end{equation}
 and vice versa.

Obviously this  operation corresponds to the projection of the canonical curve $P_\Ga \subset \PP^{g-1}$ from the node $p_e$ to the canonical
 curve $P_{\Ga'} \subset \PP^{g-2}$.
\begin{rmk} This projection is an analog of the Fano double projection for Fano threefolds.
\end{rmk}

So, we have the correspondence
\begin{equation}
C \subset \sP_{g} \times \sP_{g-1}
\end{equation}
with two projections:
\begin{equation}
p_g \colon C \to \sP_{g}
\end{equation}
with the fiber 
$$
p_g^{-1} (\Ga) = E(\Ga)
$$
and
\begin{equation}
p_{g-1} \colon C \to \sP_{g-1}
\end{equation}
with the fiber 
$$
p_{g-1}^{-1} (\Ga') = S^2(E(\Ga')).
$$

 Now  we blow down this edge in our graph. We
get  a new graph with 4-valent vertex $v_{new} = v = v'$ with the
star $S(v_{new}) = e_{1}, e_{2}, e'_{1}, e'_{2}$. There are 3
partitions of this set to pairs: the old one $(e_{1}, e_{2}) \vert
(e'_{1}, e'_{2})$ and two new ones:
\begin{equation}
(e_{1}, e'_{2} \vert e_{2}, e'_{1}) \text{ and } (e_{1}, e'_{1}
\vert e_{2}, e'_{2}).
\end{equation}
Now we can blow up the vertex  $v_{new}$ to the edge $e_{new}$
with vertices $\p e_{new} = v_{new}, v'_{new}$ and stars
\begin{enumerate}
\item $S(v_{new}) = e_{new}, e_{1}, e_{2}$ and $S(v'_{new}) = e_{new},
e'_{1}, e'_{2}$. This is our starting graph $\Ga$.
\item $S(v_{new}) = e'_{new}, e_{1}, e'_{2}$ and $S(v'_{new}) = e'_{new},
e'_{1}, e_{2}$. This is the first new graph $\Ga'$.
\item $S(v_{new}) = e''_{new}, e_{1}, e'_{1}$ and $S(v'_{new}) = e''_{new},
e'_{2}, e'_{2}$. This is the second new graph $\Ga''$.
\end{enumerate}
As a result of this construction every obtained graph
 has distinguished edge that is we have a triple of flags
\begin{equation}
(e \subset \Ga), ( e'_{new} \subset \Ga'), ( e''_{new} \subset
\Ga'').
\end{equation}
Such triple we call a {\it nest} of flags.

The correspondence    (2.14) gives a nest of graphs of genus $g-1$
\begin{equation}
(e_1, e'_1 \subset \Ga_1), (e_2, e'_2 \subset \Ga_2), ((e_3, e'_3 \subset \Ga_3))
\end{equation}
 Every nest is defined
uniquely by any flag $(e \subset) \Ga$  or $(e_i, e'_i \subset \Ga_i)$ from the triple.

 Moreover,
let $S(e)$ be the  star  of an edge $e$ that is  the union
of stars of the boundary $v, v' = \p e$
\begin{equation}
S(e) = S(v)\bigcup S(v').
\end{equation}
Then we have the canonical identification of graphs
\begin{equation}
\Ga - S(e) = \Ga' - S(e_{new}) = \Ga'' - S(e_{new}).
\end{equation}

\begin{rmk}
It is easy to see that if $e$ is a loop this construction gives
the same graph with the same loop $e$ again.
\end{rmk}

\section{Special 1-parameter deformations of large limit curves}

Now we are ready to construct 1-parameter family of
reducible curves with $2g-3$ rational components. 

For our flag $e \subset \Ga$ consider two componets $P_{v}, P_{v'}$ of the curve $P_{\Ga}$ with
 common point $p_{e}$ where $v, v' = \p e$. Remove this point and
 glue $P_{v}$ and $P_{v'}$ by a tube that is, consider the
 connected sum
 \begin{equation}
 P_{v}\#_{p(e)} P_{v'} = P_{v,v'} = S^{2}.
 \end{equation}
This is a 2-sphere with two pairs of points $(p_{e_{1}}, p_{e_{2}})$
and $(p_{e'_{1}}, p_{e'_{2}})$ where $(e, e_{1}, e_{2}) = S(v)$
and $(e, e'_{1}, e'_{2}) = S(v')$. If we fix a complex structure
on $S^{2}$ and consider the double cover
\begin{equation}
\phi \colon E \to \PP^{1}
\end{equation}
with ramification points
\begin{equation}
W = p_{e_{1}} \bigcup p_{e_{2}} \bigcup p_{e'_{1}} \bigcup
p_{e'_{1}}
\end{equation}
we obtain an elliptic curve $E$ with a point of second order
\begin{equation}
\si = p_{e_{1}} + p_{e_{2}} - p_{e'_{1}} - p_{e'_{2}} \in Pic(E)_2.
\end{equation}
So the moduli space of complex structures on $S^{2}$ with quadruple of points divided in two pairs is equal to
$\sM_{1}^{2}$ - the moduli space of smooth elliptic curves with fixed
point of order 2.

Every such complex structure $\tau \in \sM_{1}^{2}$ on $S^{2}$  with quadruple of points divided in two pairs  and the
standard complex structures on all others components glued as before define a
stable algebraic reducuble curve $P_{e \subset\Ga, \tau}$. Thus
we obtain an embedding
\begin{equation}
\psi_{e \subset \Ga} \colon \sM_{1}^{2} \to \ov{\sM_{g}}.
\end{equation}
The moduli space
\begin{equation}
\sM_{1}^{2} = \PP^{1} - (\si, \si', \si'')
\end{equation}
is the projective line without 3 points. These 3 points
correspond to 3 possibility to divide 4 points $p_{e_{1}},
p_{e_{2}}, p_{e'_{1}}, p_{e'_{2}}$ in 2 pairs that is a choice of
a point of order 2 on an elliptic curve.

It is easy to see that

\begin{prop}
\begin{enumerate}
\item The map $\psi_{e \subset \Ga}$   (3.5) can be extended to a map
\begin{equation}
\psi_{e \subset \Ga} \colon \PP^{1} \to \ov{\sM_{g}};
\end{equation}
\item 
\begin{equation}
\psi_{e \subset \Ga} (\si) = p_{e} \in  P_{\Ga}
\end{equation}
where $\si$ is given by   (3.4) and $P_{\Ga}$ is a large limit curve with
fixed node corresponding to the edge $e$,
\item
\begin{equation}
\psi_{e \subset \Ga} (\si') = P_{\Ga'};
\end{equation}
and
\item
\begin{equation}
\psi_{e \subset \Ga} (\si'') = P_{\Ga''};
\end{equation}
where the  triple
\begin{equation}
(e \subset \Ga), (e_{new} \subset \Ga'), ( e_{new} \subset \Ga'')
\end{equation}
is a nest of flags    (2.17).
\item Now we can identify the tangent direction
\begin{equation}
T\psi_{e \subset \Ga} (\PP^{1})_{P_{\Ga}} = \C_{e}
\end{equation}
from the decomposition    (2.11) with the corresponding node of the
double canonical curve. We obtain the uniquely defined  2-canonical  model.
\end{enumerate}
\end{prop}

 This identification of images of nodes
under the double canonical embedding and the directions of special deformations cancels  the projective transformation ambiguity.

\begin{rmk} If  edge $e$ is a loop with a vertex $v$ such that
the star $S(v) = e, e'$ then $P_{v}$ is a rational curve with one
double point $p_{e}$ and the smooth point $p_{e'}$. The operation
of connected summing  around double point $p_{e}$ gives a smooth
2-torus $T^{2}$ with fixed point $p_{e'}$ and the isotopy class
$a \in H_{1}(T^{2}, \Z)$ which is the class of the neck of gluing
tube. The class $a$ mod 2 gives a point of order 2 on $T^{2}$.
The space $\sM_{1}^{2}$ of complex structures on $T^{2}$  with
such additional data is $\C^{*} - 1$ that is $\PP^{1}$ without 3
points again. In this case the rational curve $\psi_{e \subset
\Ga} (\sM_{1}^{2})$ admits the compactification by the double
point corresponding to $P_{\Ga}$.
\end{rmk}

So, the Deligne-Mumford  compactification $\ov{\sM_{g}}$ contains
the configuration $\sC$ of rational curves
\begin{equation}
\{ C \} = \sC,  \text{ where every } C = \psi_{e \subset \Ga} (\PP^{1})
\end{equation}
for some flag $e \subset \Ga$. For three flags from the same nest    (2.17) the curve $C$ is the same.  We can consider this configuration of rational
curves as a reducible curve
\begin{equation}
\sC = \bigcup C
\end{equation}
that is the union of all components. It is easy to see that

\begin{prop}
\begin{enumerate}
\item for every  component $C$ of $\sC$
\begin{equation}
  C \bigcap \sP_g  = P_{\Ga} + P_{\Ga'} + P_{\Ga''}
\end{equation}
where $(\Ga, \Ga', \Ga'')$ is a nest of graphs    (2.17);
\item for a pair $C, C'$ of componets  either the intersection $C \bigcap C'$
 is empty or it is transversal and
 \begin{equation}
 C \bigcap C' \in \sP
 \end{equation}
 \item the set $S(\Ga)$ of components through every point
  $P_{\Ga} \in \sP$ is enumerated by the set $E(\Ga)$;
\end{enumerate}
\end{prop}

Now let $ETG_g$ be the set of 3-valent flags of genus $g$ with
 the projection to the set of graphs
\begin{equation}
\ga \colon ETG_g \to TG_g, \quad \ga^{-1}(\Ga) = E(\Ga)
\end{equation}
and $Com(\sC)$ be the set of components of the reducible curve
 $\sC$. Then we have 3-cover
\begin{equation}
c \colon ETG_g \to Com(\sC)
\end{equation}
with remification along flags $e \subset \Ga$, where $e$ is a loop in $\Ga$. Let $LTG_g$ be the set of such loop flags. Then
one has
\begin{equation} 
3 \cdot \vert  Com(\sC) \vert - \vert LTG_g \vert = 3(g-1) \cdot \vert TG_g \vert.
\end{equation}
Thus, roughly speaking, the cardinality of $Com(\sC)$ is equal to $(g-1)$
times the cardinality of $\sP_g$.

We saw that if an edge $e \in E(\Ga)$ is not a loop then the map
$\psi_{e \subset \Ga}$   (3.7) is an embedding and the corresponding component $C$
 of the reducible curve $\sS \subset \ov{\sM_{g}}$ is a
smooth rational curve with fixed 3 points $ ( P_{\Ga}, P_{\Ga'},
P_{\Ga''})$ such that the corresponding 3 graphs  can be equipped with flags structures forming the nest
   (2.17). Over every of such points the tangent space admits
decomposition    (2.11). 

Geometry of the other curves from the family
of curves parametrized by $C$ is very near to the geometry of
large limit curves: let $e \in E(\Ga)$, $v, v' = \p e$ and $P_{E,
\si}$ is the point of $C$ corresponding to the elliptic curve
  (3.2) with a point of order 2   (3.4). Then $P_{E, \si}$ has 2g-4
 old components
 \begin{equation}
 \bigcup_{v'' \neq v, v'} P_{v''}
 \end{equation}
with triples of points and 3g-1 nodes $p_{e'}, e \neq e'$ and one
new component $P_{v, v'}$  with the quadruple of points (see
  (3.1). Then

\begin{enumerate} \item
the canonical class $K_{P_{E, \si}}$ is a line bundle: a
restriction of it to every component $P_{v''}, v'' \neq v, v' $
is the sheaf of meromorphic differentials $\om$ with simple poles
at $p_{e}, p_{e'}, p_{e''}$  where $\{e,e', e''\} = S(v)$;
\item the
restriction of the canonical class to the component $P_{v, v'}$
is the sheaf of meromorphic differentials $\om$ with simple poles
at $p_{e_{1}}, p_{e_{2}}, p_{e'_{1}}, p_{e''_{2}}$ (see   (3.1) and   (3.3)).
\item
Thus 
\begin{equation}
c (K_{P_{E, \si}}) = (2,1,...,1) \in NS_{P_{E, \si}}
\end{equation}
where the first coordinate corresponds to $P_{v, v'}$.
\item Again
every holomorphic section $s$ of the canonical class  is a
collection of meromorphic differentials $\{ \om_{v''}\}$ on the
components $P_{v''}$ with poles at $p_{e}, p_{e'}, p_{e''}$ where
$\{e,e', e''\} = S(v)$ and a meromorphic differential $\om_{v,
v'}$ on the component $P_{v, v'}$ with poles at the quadruple
with the same constraints    (2.3) and    (2.4).
\item
The canonical map defined by the complete linear system $\vert K_{P_{E,
\si }}\vert$
\begin{equation}
\phi_{K} \colon P_{E, \si} \to \PP^{g-1}
\end{equation}
sends $P_{v''}$ to a configuration of lines and $P_{v, v'}$ to a
conic  in $\PP^{g-1}$.
\item Again the dimension of the space of quadratic differentials
on $P_{E, \si}$ is equal to 3g-3 and this curve is an orbifold
smooth point of $\ov{\sM_{g}}$.
\item
The double  canonical map of $P_{E, \si}$  given  by the complete
linear system $\vert 2K_{\Ga}\vert$ is an embedding
\begin{equation}
\phi_{2K_{P_{E, \si}}} \colon P_{E, \si} \to \PP^{3g-4} = \PP
T\ov{\sM_{g}}_{P_{E, \si}}.
\end{equation}
(see (2.11)).
\item The images  of nodes
\begin{equation}
\{ \phi_{2K_{P_{E, \si}}}(p_{e'}) \}, \quad e \neq e' \in
E(\Ga)
\end{equation}
define the decomposition of the restriction of the tangent bundle
\begin{equation}
T\ov{\sM_{g}} \vert_{C} = TC \bigoplus (\bigoplus_{e \neq e' \in
E(\Ga)} L_{e'})
\end{equation}
where the fiber of the line bundle $L_{e'}$ over a point is the
component of the decomposition    (2.11).
\item Thus every line bundle $L_{e'}$ from the previous decomposition
is the tautological line bundle
\begin{equation}
L_{e'} = \Oh_{C} (-1),
\end{equation}
\item hence the previous decomposition is
\begin{equation}
T\ov{\sM_{g}} \vert_{C} = \Oh_{C}(2)\bigoplus (3g-4) \Oh_{C}(-1).
\end{equation}
\item The restriction to $C$ of the canonical class of $\ov{\sM_{g}}$
is
\begin{equation}
K_{\ov{\sM_{g}}} \vert_C = \Oh_{C}(3(g-2)).
\end{equation}
\end{enumerate}

For a singular curve the moduli spaces of  vector bundles are not compact. These moduli spaces admit the compactifications by torsion free sheaves. These sheaves are not local free  only over
nodes and the theory of the compactification is very close to
the theory for algebraic surfaces. All details of this theory
can be found in forthcoming  Artamkin paper \cite{A}. All constructions are very close to the complex gauge theory on
smooth compact Riemann surfaces. We will see this in the following section where we apply  the main constructions of the complex gauge theory from section 1 to ll-curves.

\section{The complex gauge theory on ll-curves}

A vector bundle $E$ on $P_\Ga$ is called topologicaly trivial if
the restrictions $E \vert_{P_v}$ are trivial for all $v \in V(\Ga)$.
We begin with the   description of the moduli spaces of topologicaly trivial 
\begin{enumerate}
\item line bundles $Pic_0 (P_\Ga)$,
\item rk 2 semi-stable bundles $M_{vb}^{ss}$.
\end{enumerate}
Since each  component $P_{v}$ is  projective line the restriction of any topologicaly trivial line bundle equals
\begin{equation}
L \vert_{P_{v}} = \Oh_{P_{v}}.
\end{equation}

To describe a line bundle on $P_\Ga$ consider a collection of any trivializations of $L$ on all components $P_v$ of $P_\Ga$ and denote the line bundle with such addition structure by $L_0$.

Now to get a line bundle $L$ on $P_{\Ga}$ we have to concord every
 pair of line bundles $\Oh_{P_{v}}, \Oh_{P_{v'}}$
at common point $p_{e}$ if $v, v' = \p e$. Under our trivializations such concordance is
given by a multiplicative constant $ a (\vec e) \in \C^{*}$ of oriented edge
such that under the orientation reversing involution $i_e$
\begin{equation}
a(i_e(\vec e)) = a(\vec e)^{-1}.
\end{equation}
 Thus $L_0$ defines a map
\begin{equation}
a \colon \vec E(\Ga) \to \C^*.
\end{equation}
subjecting (4.2).

The changing of trivialization is given by a function
\begin{equation}
\tilde{\la} \colon V(\Ga) \to \C^*
\end{equation}
 which acts on the functions $a$  by the formula
\begin{equation}
\tilde{\la}(a(\vec e)) = \tilde{\la}(v_s) \cdot a(\vec e) \cdot \tilde{\la}(v_t)^{-1}
\end{equation}
where $\p \vec e = v_s \bigcup v_t$ and $v_s$ is the source of arrow and $v_t$ is the target.

Geometrically the construction of a concordance $a(\vec e)$
at a point $p_e$ for ll-curve plays the role of the period (monodromy) $\rho([C_e])$ of a flat connection $[\rho] \in \CLRep(\pi_1(\Si_\Ga), \C^*)$ along the cycle $\vec C_e$ for non-singular curve $\Si_\Ga$ from the previous section. 

Hence if
 $\sA_{\C}$ is the space of trivialized topologicaly trivial line bundles that is the space of functions $a$ (4.2) subjecting (4.3)
 and $\sG_{\C}$ be the group of changing of trivializations (4.4). Then the moduli space of line bundles on 
$P_\Ga$
\begin{equation}
Pic_0(P_\Ga) = Hom(ker \quad r, \C^*) = (\C^*)^g
\end{equation}
 is the abelian Schottky space.

For rk 2 bundle $E$ on $P_\Ga$ we have the same type description:
consider a collection of  trivializations  of restrictions $E \vert_{P_v} = \Oh \oplus \Oh$ on all components  of $P_\Ga$ and denote the  bundle with such addition structure by $E_0$.

Now  we have to glue every
 pair of  bundles $\Oh_{P_{v}} \oplus \Oh_{P_{v}}, \Oh_{P_{v'}} \oplus \Oh_{P_{v'}}$
at common point $p_{e}$ if $v, v' = \p e$. Such gluing is given by a function
\begin{equation}
a \colon \vec E(\Ga) \to SL(2, \C)
\end{equation}
subjection to the equation
\begin{equation}
a(i_e(\vec e)) = a(\vec e)^{-1}.
\end{equation}
The changing of the trivialization is given by the function
\begin{equation}
\tilde{g} \colon V(\Ga) \to SL(2, \C)
\end{equation}
 which acts on the functions $a$  by the formula
\begin{equation}
\tilde{g}(a(\vec e)) = \tilde{u}(v_s) \cdot a(\vec e) \cdot \tilde{g}(v_t)^{-1}
\end{equation}
where $\p \vec e = v_s \bigcup v_t$.

We would like to emphisize again that geometrically the construction of a concordance $a(\vec e)$
at a point $p_e$ for ll-curve plays the role of the period (monodromy) $\rho([C_e])$ of a flat connection $[\rho] \in \CLRep(\pi_1(\Si_\Ga), \SL(2, \C))$ along the cycle $\vec C_e$ for non-singular curve $\Si_\Ga$ from the previous section.

Again if    $\sA_{ \C}$ is the space  of trivialized topologicaly trivial rk 2 bundles with  the group of trivializations $\sG_{ \C}$ acting by the formula (4.10)
then the moduli space of semi-stable topologicaly trivial
vector bundles on $P_\Ga$ is the quotient
\begin{equation}
M^{ss}_{vb}(P_\Ga) = \sA_{\C} / \sG_{\C}.
\end{equation}
Again it is easy to see that 
\begin{equation}
M^{ss}_{vb}(P_\Ga) = \CLRep(ker \quad r, SL(2, \C)) = S_g
\end{equation}
(see     (1.31)) is the  Schottky space of genus $g$.

Now we would like to describe the space  of holomorphic flat connections on vector bundles over this
 moduli space.

To do this let us draw the trinion decomposition corresponding
to $\Si_\Ga$     (1.1)-    (1.2) again. We can construct the space of
holomorphic flat bundles 
 gluing the spaces of $\SL(2, \C)$-flat connections on the every component  $P_v - \bigcup_{e \in S(v)} p_e$ that is on $\PP^1$ without 3 points.

Such spaces of $\SL(2, \C)$-flat connections are given as the classes of representations spaces
\begin{equation}
A^{na}(P_v - \bigcup_{e \in S(v)} p_e) =  CLRep(\pi_1(P_v - \bigcup_{e \in S(v)} p_e), \SL(2, \C))
\end{equation}
where $\pi_1(\tilde{v}) = F_2$ is the free group with 2 generators. Thus we can identify the spaces
\begin{equation}
A^{na}(P_v - \bigcup_{e \in S(v)} p_e) = A^{na}(\tilde{v}) = S_2
\end{equation}
(see     (1.30)) where the last one is the Schottki space of genus 2.

If two vertices $v$ and $v'$ are joined by the
 edge $e$ that is if the components $P_v$ and $P_{v'}$ intersect at $p_e$ then we can glue $A^{na}(P_v - \bigcup_{e \in S(v)} p_e) $ and $A^{na}(P_{v'} - \bigcup_{e \in S(v')} p_e) $ along this intersection point
 $p_e$ 
 just in the same way as we do this for the smooth case     (1.32)-    (1.35). Again for a semi-simple element
$m \in \SL(2, \C)$ the stabilizer $Z(m) = \C^*$.  The result of such gluing is
\begin{equation}  
A^{na}(P_v - \bigcup_{e \in S(v)} p_e) * A^{na}(P_{v'} - \bigcup_{e \in S(v')} p_e)) = S_3
\end{equation}
if $P_v \bigcap P_{v'} = p_e$.
Gluing all components we get our ll-curve $P_\Ga$ and gluing all spaces of flat bundles $A^{na}(P_v - \bigcup_{e \in S(v)} p_e))$ we get the same space $ \CLRep(\pi_1(\Si_\Ga), \SL(2, \C))
$ of flat bundles on $P_\Ga$ as for the smooth case:
\begin{equation}
A^{na}(P_\Ga) = \CLRep(\pi_1(\Si_\Ga), \SL(2, \C))
 = A^{na}(\Si_\Ga).
\end{equation}
The direct interpretation gives the following coincidence:

\begin{prop} 
\begin{enumerate}
\item The forgetful map     (1.25) for ll-curve $P_\Ga$
\begin{equation}
f \colon A^{na}(P_\Ga) =  \CLRep(\pi_1(\Si_\Ga), \SL(2, \C)) \to M^{ss}_{vb}(P_\Ga) = 
\end{equation}
$$
= \CLRep(ker \quad r, \SL(2, \C)) = S_g
$$
is the natural map induced by the restriction of representations to the kernel $ker \quad r \subset \pi_1(\Si_\Ga)$ (see     (1.7)-    (1.9)).
\item This forgetful map admits a holomorphic section
\begin{equation}
s \colon  \CLRep(ker \quad r, \SL(2, \C)) = S_g \hookrightarrow  CLRep(\pi_1(\Si_\Ga), \SL(2, \C)) 
\end{equation}
induced by  epimorphism $r$     (1.7) and  isomorphism
$Int$     (1.20). 
\item Thus the affine bundle (4.17) is the Higgs fields vector bundle (or the Hitchin bundle from \cite{H}), which we may consider as the cotangent   bundle to the singular variety $ M^{ss}_{vb}(P_\Ga) = S_g$.
\item For every complex structure $I \in \tau_g$ (see the text between formula     (1.20) and     (1.21)) we have the rational map
\begin{equation}
F_{(\Ga, I)} \colon  M^{ss}_{vb}(P_\Ga) = S_g \to M^{ss}(\Si_I)
\end{equation}
which is the composition of  section $s$ (4.18) and  forgetful map $f$     (1.25). 
\end{enumerate}
\end{prop}

Again  consider the space of representations $Rep(\pi_1(\Si), \SL(2, \C))$ ( without the diagonal adjoint factorization) (1.37) with the projection $p$
(1.38).
The geometrical  meaning of this procedure is  almost the same as before: 
 we fixed a component $P_v$ and a trivialization 
$E \vert_{P_v} = \C^2$  of the vector bundle $E$ over this component. The moduli space of vector bundles with such additional structure 
\begin{equation}
\widetilde{M^{ss}_{vb}(P_\Ga)} = (\SL(2,\C))^g
\end{equation} 
is a non singular algebraic variety with
 the structure of  principal $PGL(2, \C)$-bundle over $S_g$: 
\begin{equation}
\phi   (\SL(2,\C))^g  \to S_g
\end{equation}
where the group $\SL(2, \C)$ modulo $\pm 1$ acts on  $E  \vert_{P_v} = \C^2$ as (2 X 2)-matrices with determinant 1.  

The moduli space of flat vector bundles on $P_\Ga$ with such additional structure is $Rep(\pi_1(\Si), \SL(2, \C))$ with the structure $p$ (1.38) of a principal $PGL(2, \C)$-bundle. The moduli space of flat bundles  $Rep(\pi_1(\Si), \SL(2, \C))$ can be described as the preimage of 1 for the algebraic map 
\begin{equation}
com \colon (\SL(2,\C))^{2g} = \prod \SL(2, \C)_{a_i}
\times  \prod \SL(2, \C)_{b_i} \to \SL(2, \C)
\end{equation}
$$
com ( g(a_1), ... , g(a_g),  g(b_1), ... , g(b_g)) =
\prod_{i=1}^g [g(a_i), g(b_i)] 
$$
(see (1.9)):
$$
com^{-1} (1) = Rep(\pi_1(\Si), \SL(2, \C))
$$
which is an algebraic variety obviously.

Again  regular functions (1.40)-(1.41) defines the divisors (1.42) with the intersection (1.43).

But now  the affine bundle     (4.17) over $M^{ss}_{vb}(P_\Ga)$ is a vector bundle because it admits the holomorphic section (4.18). It has as the fiber over  a
 point $E \in M^{ss}_{vb}(P_\Ga)$ the space
\begin{equation} 
f^{-1}(E) = H^0(P_\Ga, ad E \times K_\Ga) =\C^{3g-3}_E
\end{equation}
which is  a {\it packet} of classes of representations of $\pi_1(\Si_\Ga)$ (see (1.2), (1.6)).

Again consider the restriction of the principal $\SL(2, \C)$-bundle
(1.38) to this affine space
\begin{equation}
p \colon  p^{-1} ( H^0(P_\Ga, ad E \times K_\Ga) ) \to  H^0(P_\Ga, ad E \times K_\Ga).
\end{equation}
Restricting the functions (1.40) to this affine variety
 we get the complete intersection
\begin{equation}
 i(S_g) \cap  ( H^0(P_\Ga, ad E \times K_\Ga) ) =  E
\end{equation}
as a point of the moduli space $M^{ss}_{vb}(P_\Ga) = S_g$. That is for every vector bundle $E$ on $P_\Ga$
\begin{enumerate}
\item {\it this complete intersection is non empty}, and
\item {\it there exists unique class of the Schottki representation.} 
\end{enumerate}

From this we get immideately our final result:

\begin{thm} On  general curve  general stable vector bundle admits the Schottki representation.
\end{thm} 
\begin{rmk} Of course all our constructions are valid for
 vector bundle of any rank. For simplicity we  worked  with rk 2.
\end{rmk}

In this paper we  uses only direct "classical" arguments. But obviously these results can be refined by using the Higgs bundles technique for ll-curves. In particular 
 we know (see \cite{H}) that the Higgs pairs space $ T^*M^{ss}(P_\Ga)$ admits a partial compactification
$\ov{T^*M^{ss}(P_\Ga)}$ completing the fibers of the holomorphic moment map.
This partial compactification defines a compactification
\begin{equation} 
\ov{S_g} \subset  \ov{T^*M^{ss}(P_\Ga)}.
\end{equation}
of the moduli space of vector bundles on ll-curves by torsion free sheaves. Now  we postpone this beautifull subject to the next  paper of the serie refining \cite{T3}.

\end{document}